\documentclass[12pt]{amsart}
\pdfoutput=1
\usepackage{amsmath, amssymb}
\usepackage{amsfonts}
\usepackage{mathrsfs}
\usepackage[arrow,matrix,curve,cmtip,ps]{xy}
\usepackage{graphicx}
\usepackage{amsthm}
\usepackage{float}
\usepackage{amsthm}
\usepackage[utf8]{inputenc}
\usepackage[T1]{fontenc}
\usepackage{mathtools}

\allowdisplaybreaks

\newtheorem{theorem}{Theorem}[section]

\newtheorem{corollary}[theorem]{Corollary}

\newtheorem*{theorem*}{Theorem}
\theoremstyle{remark}
\newtheorem{remark}[theorem]{Remark}

\newtheorem{question}{Question}

\numberwithin{equation}{section}

\begin{document}
\title[Lower bounds for $k$-systems]{ Constructing large $k$-systems on surfaces}

\author{Tarik Aougab}

\address{Department of Mathematics \\ Yale University \\ 10 Hillhouse Avenue, New Haven, CT 06510 \\ USA}
\email{tarik.aougab@yale.edu}

\date{\today}


\begin{abstract}

Let $S_{g}$ denote the genus $g$ closed orientable surface. For $k\in \mathbb{N}$, a \textit{k-system} is a collection of pairwise non-homotopic simple closed curves such that no two intersect more than $k$ times. Juvan-Malni\v{c}-Mohar \cite{Ju-Mal-Mo} showed that there exists a $k$-system on $S_{g}$ whose size is on the order of $g^{k/4}$. For each $k\geq 2$, We construct a $k$-system on $S_{g}$ with on the order of $g^{\lfloor (k+1)/2 \rfloor +1}$ elements. The $k$-systems we construct behave well with respect to subsurface inclusion, analogously to how a pants decomposition contains pants 
decompositions of lower complexity subsurfaces. 

\end{abstract}

\maketitle

\section{Introduction}
Let $S_{g,p}$ denote the compact orientable surface of genus $g$ with $p$ boundary components. A \textit{k-system} is a collection of essential, pairwise non-homotopic simple closed curves $\left\{\gamma_{1},...,\gamma_{n}\right\}$ on $S_{g,p}$ such that no two curves in the collection intersect more than $k$ times.

Let $N(k,g,p)$ denote the maximum cardinality of a $k$-system on $S_{g,p}$, and let $N(k,g):= N(k,g,0)$. Juvan-Malni\v{c}-Mohar \cite{Ju-Mal-Mo} first showed that for any pair $(k,g), N(k,g)<\infty$. Furthermore, they produce lower bounds which grow asymptotically like $g^{k/4}$. Concretely, they show:

\begin{theorem} \cite{Ju-Mal-Mo} Given $k \in \mathbb{N}$, for sufficiently large genus $g$, there exists a $k$-system on $S_{g}$ of size at least ${n \choose \lfloor k/2 \rfloor}$, where 

\[ n = \sqrt{25+ 48(g-1)}-5/2.\]

\end{theorem}

The main focus of this article is to improve these lower bounds for all $k\geq 2$ by constructing large $k$-systems. Specifically, we show:

\begin{theorem}  Given $k\in \mathbb{N}$, $k\geq 2$, there exists a $k$-system $\Omega(k,g)$ on $S_{g}$ such that 

\[ |\Omega(k,g)| \geq {g \choose 1+\lfloor\frac{k}{2} \rfloor} ,\]

and 
\[\Omega(k,g)= \Theta(g^{ \lfloor (k+1)/2 \rfloor +1}).\]

\end{theorem}

In the statement of Theorem $1.2$, $\Theta()$ denotes asymptotic growth rate:

\[ f= \Theta(h) \Leftrightarrow 0< \lim_{n\rightarrow \infty} \frac{f(n)}{h(n)} <\infty. \] 

Note that our lower bound is not asymptotic, in the sense that it does not require $g$ to be sufficiently large with respect to $k$.

\begin{remark} We do not expect that these lower bounds are sharp. Indeed, the lower bound for $k=2$ grows like $g^{2}$, and Malestein-Rivin-Theran \cite{Mal-Riv-Ther}, and independently Constantin \cite{Con} have constructed $1$-systems on $S_{g}$ with quadratically many elements. 

\end{remark}

\begin{remark} Przytycki \cite{Prz} has recently shown an upper bound for the maximum size of a $1$-system which grows cubically in $|\chi(S)|$, for $\chi$ the Euler characteristic. In the same paper, for each $k>1$, he obtains an upper bound for the maximum size of a $k$-system which grows like $|\chi|^{k^{2}+k+1}$. Alternatively, Juvan-Malni\v{c}-Mohar \cite{Ju-Mal-Mo} provide an upper bound for the size of a $k$-system on $S_{g}$ which grows like 
\[ 2^{\chi}\cdot [\chi^{2}k]^{\chi k}.\]
For fixed $k$, as $\chi \rightarrow \infty$, Przytycki's upper bound grows slower than Juvan-Malni\v{c}-Mohar's; however for fixed $\chi$, as $k\rightarrow \infty$, the upper bound of  Juvan-Malni\v{c}-Mohar grows slower than Przytycki's. In either case, there is still a very large gap between the size of the $k$-systems constructed here and the best known upper bounds.

\end{remark}

A maximal $0$-system is simply a pants decomposition of $S_{g,p}$. Moreover, if $c$ is a non-separating simple closed curve in a pants decomposition $\mathcal{P}$, then cutting along $c$, gluing in disks along the resulting two boundary components, and deleting any of the remaining curves in $\mathcal{P}$ which have become homotopically trivial, or homotopically redundant (i.e., there may exist distinct elements of $\mathcal{P}$ which become homotopic during this process), one obtains a subcollection $\mathcal{P}'$ which is a pants decomposition on a lower genus surface. 

The $k$-system $\Omega(k,g)$ on $S_{g}$ that we construct satisfies an analogous property in the sense that it ``contains'' $k$-systems on lower complexity subsurfaces; moreover, it also contains large $j$-systems for any $j<k$:

\begin{theorem} The $k$-systems $\Omega(k,g)$ satisfy the following properties:

\begin{enumerate}

\item $\Omega(k,g)$ contains a non-separating simple closed curve $c_{g}$ such that cutting along $c_{g}$, gluing in disks along the resulting boundary components, and deleting any element of $\Omega(k,g)$ which intersects $c_{g}$ essentially yields the collection $\Omega(k,g-1)$ on $S_{g-1}$;

\item $\Omega(k,g)$ contains a copy of $\Omega(k-1,g)$ as a subcollection.

\end{enumerate}

\end{theorem}

Let $ \mathcal{F}:= \left\{\Lambda(k,g)\right\}_{k,g}$ be a family of curve systems such that $\Lambda(k,g)$ is a $k$-system on $S_{g}$. We say that $\mathcal{F}$ satisfies property $\mathcal{I}$ (for ``inclusion'') if it satisfies the conclusions of Theorem $1.5$. As mentioned in Remark $1.3$, we do not expect our lower bounds to be best possible, but it is another interesting question to restrict attention to those families of curve systems satisfying property $\mathcal{I}$:

\begin{question} What is the maximum growth rate (in both $k$ and $g$) of a family $\mathcal{F} = \left\{\Lambda(k,g)\right\}$ satisfying property $\mathcal{I}$?

\end{question}

Furthermore, our $k$-systems have the property that for each $g$, there exists a simple closed curve $\eta$ on $S_{g}$ disjoint from $\Omega(k,g)$. As a corollary, we obtain the same lower bounds for $S_{g-1,2}$:

\begin{corollary}

$$N(k,g-1,2) \geq {g-1 \choose 1+\lfloor\frac{k}{2} \rfloor} .$$

\end{corollary}

Recall that the \textit{curve graph} of $S_{g,p}$, denoted $\mathcal{C}(S_{g,p})$, is a locally infinite, infinite diameter $\delta$-hyperbolic graph \cite{Mas-Mur} whose vertex set corresponds to the set of all isotopy classes of essential simple closed curves on $S_{g,p}$, and whose edges correspond to pairs of curves that can be realized disjointly on $S_{g,p}$. 

The fact that for each $k,g$, there exists an essential simple closed curve in the complement of $\Omega(k,g)$ implies that our $k$-systems project to diameter $2$-subsets of the corresponding curve graph.

This motivates the following question:

\begin{question}
Fix $k\in \mathbb{N}$, and let $f_{k}(g):\mathbb{N} \rightarrow \mathbb{N}$ be a function such that $\lim_{g \rightarrow \infty}f_{k}(g)=\infty$.  What is the maximum growth rate (as a function of $g$) of a family of $k$-systems $\left\{\Sigma(k,g)\right\}_{g=1}^{\infty}$, such that for each $g$, $\Sigma(k,g)$ projects to a subset of $\mathcal{C}(S_{g})$ of diameter at least $f_{k}(g)$? 

\end{question}

\begin{remark} Given $\lambda \in (0,1)$, for all $g$ sufficiently large, if a pair of curves $\alpha,\beta$ are distance at least $n$ apart in $\mathcal{C}(S_{g})$, they must intersect at least $\lceil g^{\lambda(n-2)} \rceil$ times (see \cite{Aoug}). Therefore, $f_{k}(g)$ should be chosen to have growth at most logarithmic with base $g$.

\end{remark}

\textbf{Organization of paper.} In section $2$, we introduce some basic terminology for curves on surfaces. In section $3$, we present a new method for constructing $1$-systems with quadratically many elements; these $1$-systems will serve as a sort of backbone for the $k$-systems constructed in later sections. In section $4$, we complete the general construction. 

\textbf{Acknowledgements} The author would like to thank Yair Minsky and Igor Rivin for their time and for helpful conversations during this project. He also thanks Kyle Luh and Daniel Montealegre for suggesting the use of Pascal's identity at the end of Section $4$. The author was partially supported by NSF grants DMS 1005973 and 1311844.

\section{Terminology}

\subsection{curves on surfaces.}
A curve $\gamma$ on $S_{g}$ is \textit{essential} if it is homotopically non-trivial. A \textit{multi-curve} is a collection of pairwise non-homotopic and pairwise disjoint simple closed curves. Given two homotopic curves $\gamma, \gamma'$, we write $\gamma \sim \gamma'$ for the homotopy relation. Given two homotopy classes of curves $[\alpha], [\beta]$, the \textit{geometric intersection number}, denoted $i([\alpha], [\beta])$ is simply the minimum set theoretic intersection, taken over all representatives in the homotopy classes of $\alpha$ and $\beta$:

\[ i([\alpha], [\beta]) = \min_{x \sim \alpha} |x \cap \beta|.\]

As is customary, we write $i(\alpha, \beta)$ to mean $i([\alpha], [\beta])$. If $\alpha$ is a simple closed curve on $S_{g}$, a \textit{regular neighborhood} of $\alpha$ is an embedded annulus $A$ containing $\alpha$ and which deformation retracts to $\alpha$. If $\gamma$ is an embedded arc on $S_{g}$, by a regular neighborhood of $\gamma$, we mean the image of a homeomorphic embedding $\phi: [0,1] \times [0,1] \hookrightarrow S_{g}$ such that $\phi( \left\{1/2\right\} \times [0,1] )= \gamma$.

\section{Constructing $1$-systems}
A quadratic lower bound for the maximum size of a $1$-system has been found by Malestein-Rivin-Theran \cite{Mal-Riv-Ther} and also by Constantin \cite{Con}. In this section, we construct a quadratically growing sequence of $1$-systems using a different method. These $1$-systems  will serve as a ``scaffold'' for the $k$-systems $\Omega(k,g)$ in the next section.

We begin by constructing a certain realization of $S_{g}$ that will be convenient for displaying the desired $1$-system, $\Omega(1,g)$. Beginning with $S_{1}$, recall that free homotopy classes of essential simple closed curves on $S_{1}$ are in correspondence with pairs of coprime integers. Let $\alpha_{1}$ be an arc on $S_{1}$ which runs parallel to a portion of the $(1,0)$ curve. We obtain our desired realization of $S_{2}$ by first excising a pair of small open disks $D^{(1)}_{1}, D^{(1)}_{2}$ from $S_{1}$, located within a small regular neighborhood $N_{1}$ of $\alpha_{1}$ near the endpoints of $\alpha_{1}$, and not separated by $\alpha_{1}$ within $N_{1}$. 

\begin{figure}[H]
\centering
	\includegraphics[width=2.5in]{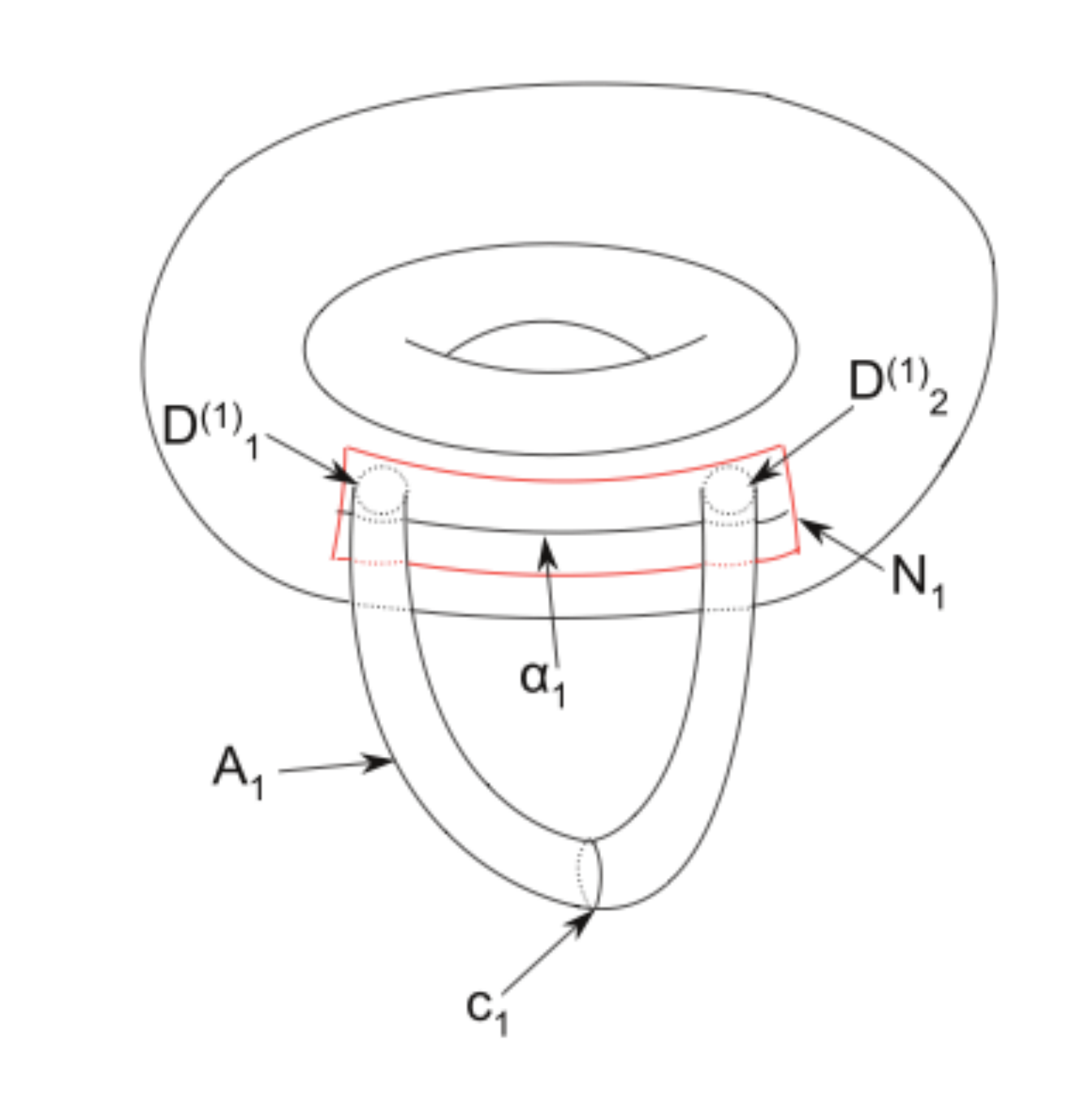}
\caption{ $\alpha_{1}$ is an arc which runs parallel to the $(1,0)$ curve on $S_{1}$, and $N_{1}$ is enclosed in a rectangle. We excise a pair of disks on the same side of $\alpha_{1}$ within $N_{1}$, and glue on an annulus $A_{1}$ with core curve $c_{1}$.}
\end{figure}

 We then glue on an annulus $A_{1}$ along the resulting two boundary components, yielding $S_{2}$; let $c_{1}$ denote the core curve of $A_{1}$. There is a simple closed curve $d_{1}$ containing $\alpha_{1}$ as a sub-arc and which intersects $c_{1}$ once, as in the figure.

\begin{figure}[H]
\centering
	\includegraphics[width=2.5in]{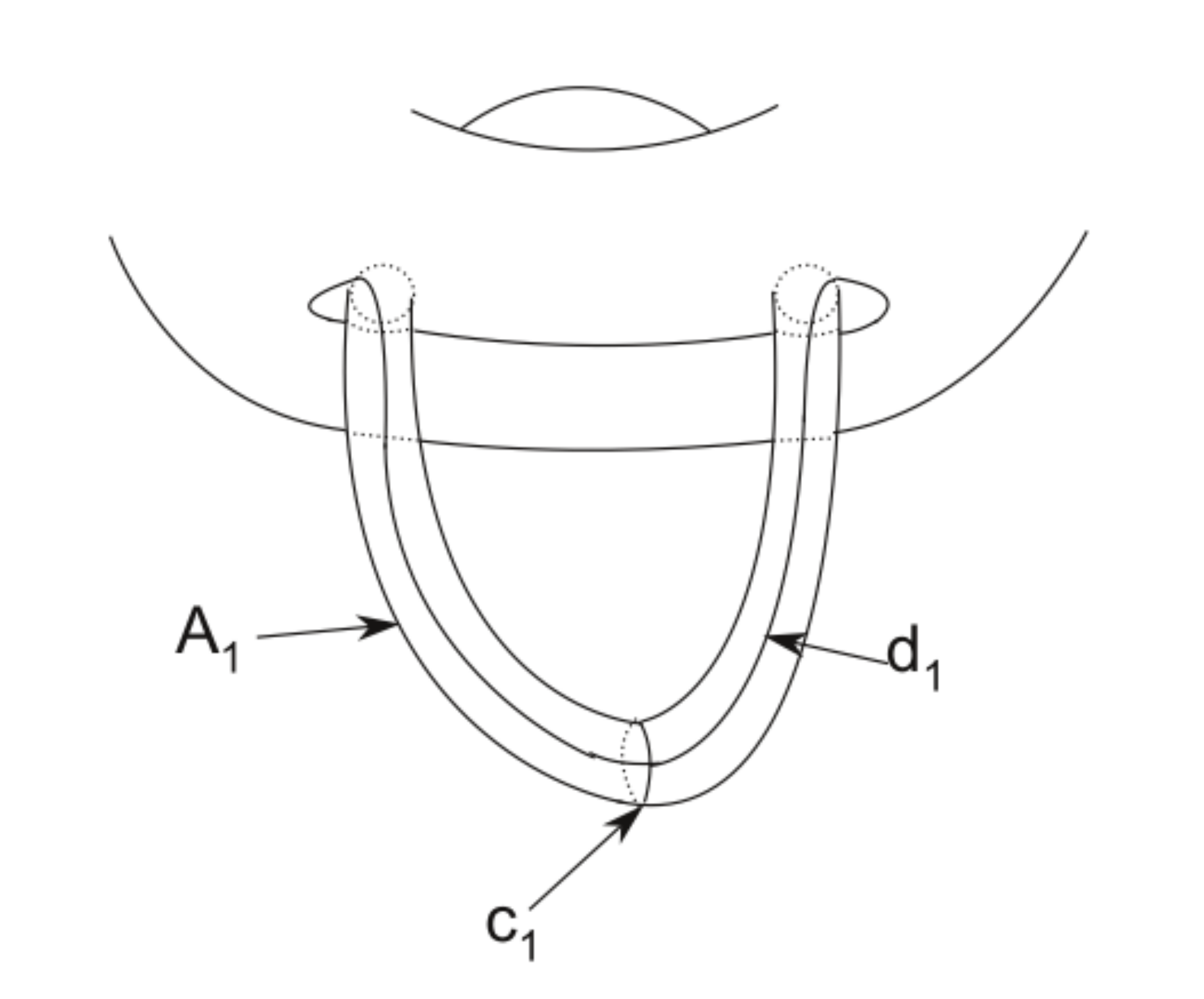}
\caption{ $d_{1}$ contains $\alpha_{1}$ as a sub-arc, and intersects $c_{1}$ once.}
\end{figure}

$d_{1}$ contains a sub-arc, $\alpha_{2}$ which is an extension of $\alpha_{1}$ and which intersects $c_{1}$ once. Consider a small regular neighborhood $N_{2}$ of $\alpha_{2}$, satisfying the following property: 

Let $N'_{2} \subset N_{2}$ denote the subset of $N_{2}$ which is a regular neighborhood of $\alpha_{1}$. Then $\partial D^{(1)}_{i}, i=1,2$ are contained in $N_{1} \setminus N'_{2}$. 

Then we obtain $S_{3}$ by excising small disks $D^{(1)}_{2}, D^{(2)}_{2}$ within $N_{2}$, not separated within $N_{2}$ by $\alpha_{2}$, and on the same side of  $\alpha_{1}$ within $N_{1}$ as $\partial D^{(1)}_{i}, i=1,2$, and gluing on an annulus $A_{2}$ along the resulting boundary components.  

\begin{figure}[H]
\centering
	\includegraphics[width=2.5in]{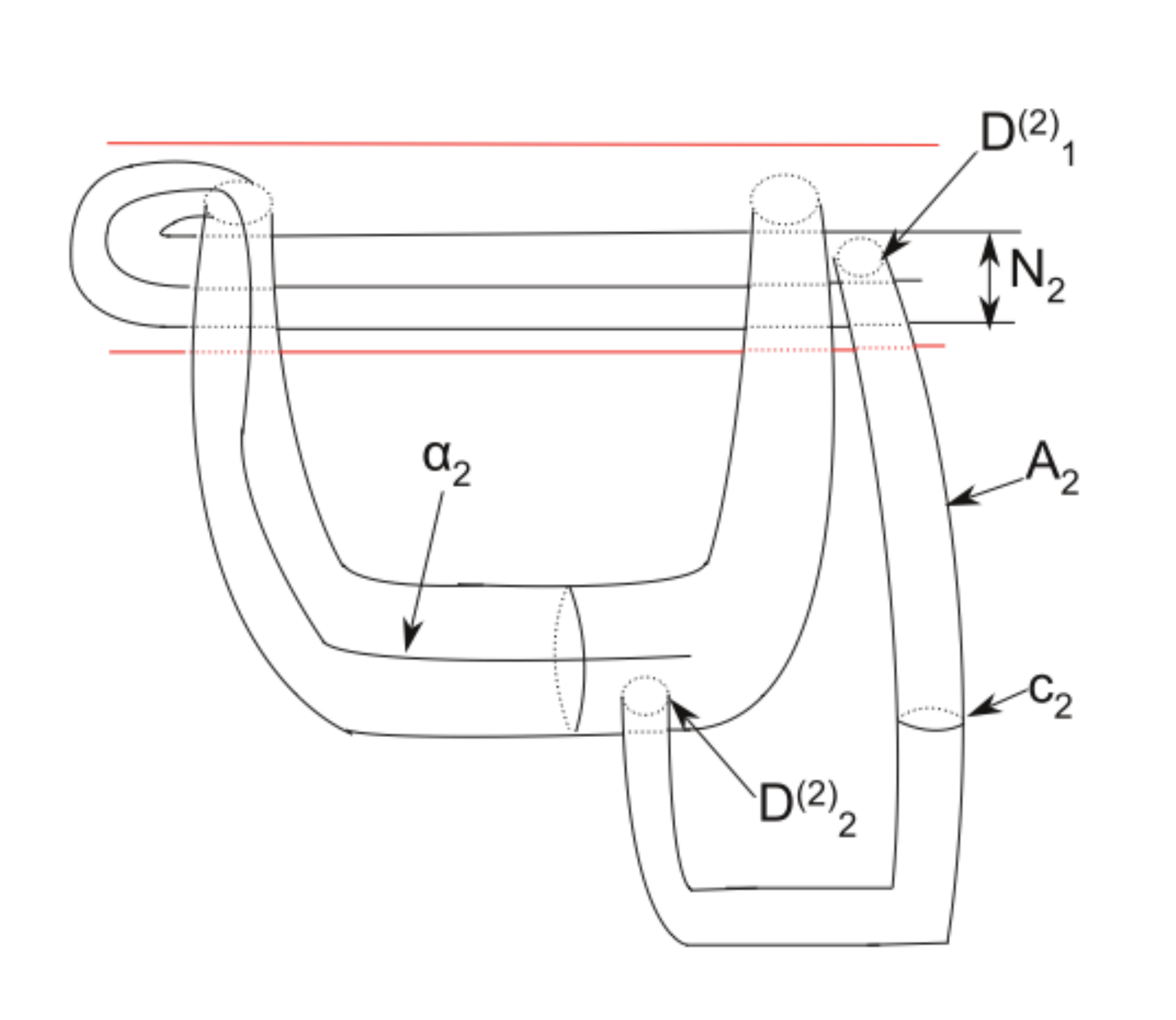}
\caption{ $N_{2}$ is closer to $\alpha_{1}$ than $N_{1}$}
\end{figure}

Note that there is a simple closed curve $d_{2}$ containing $\alpha_{2}$ as a sub-arc, and such that $i(d_{2}, c_{2})=i(d_{2},c_{1})=1$ and $i(d_{2},d_{1})=0$.  

\begin{figure}[H]
\centering
	\includegraphics[width=2.5in]{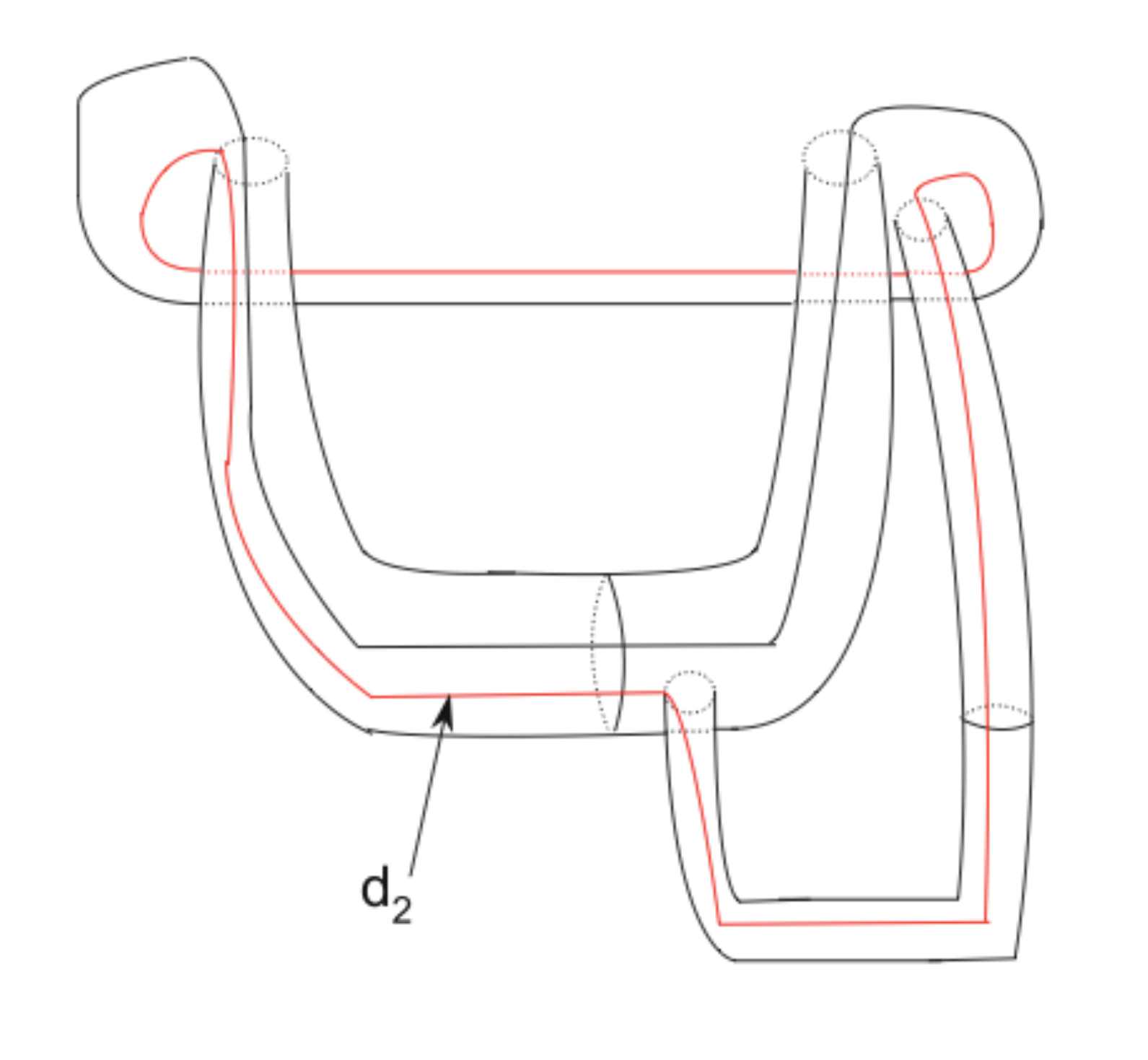}
\caption{ $d_{1}$ and $d_{2}$ }
\end{figure}

We continue inductively; on $S_{g-1}$, there is a sequence of nested arcs $\left\{\alpha_{1},....,\alpha_{g-1}\right\}$, and a collection $\left\{c_{1},...,c_{g-2}\right\}$ of pairwise disjoint simple closed curves such that $\alpha_{k}$ intersects $c_{j}$ if and only if $k \geq j-1$. There is furthermore a sequence of annuli $\left\{A_{1},..., A_{g-2}\right\}$ such that $c_{k}$ is the core curve of $A_{k}$, as well as a sequence of simply connected regions $N_{1},...,N_{g-1}$ such that $N_{k}$ is a regular neighborhood of $\alpha_{k}$. For each $k$, let $N'_{k} \subset N_{k}$ denote the subset of $N_{k}$ which is a regular neighborhood of $\alpha_{k-1}$; then for each $k$,

\[ N'_{k} \subset N_{k-1}.\]

For each $k \leq g-2$, the boundary components of $A_{k}$ are contained in $N_{k} \setminus N'_{k+1}$, and not separated by $\alpha_{k}$ within $N_{k}$.

\begin{figure}[H]
\centering
	\includegraphics[width=2.5in]{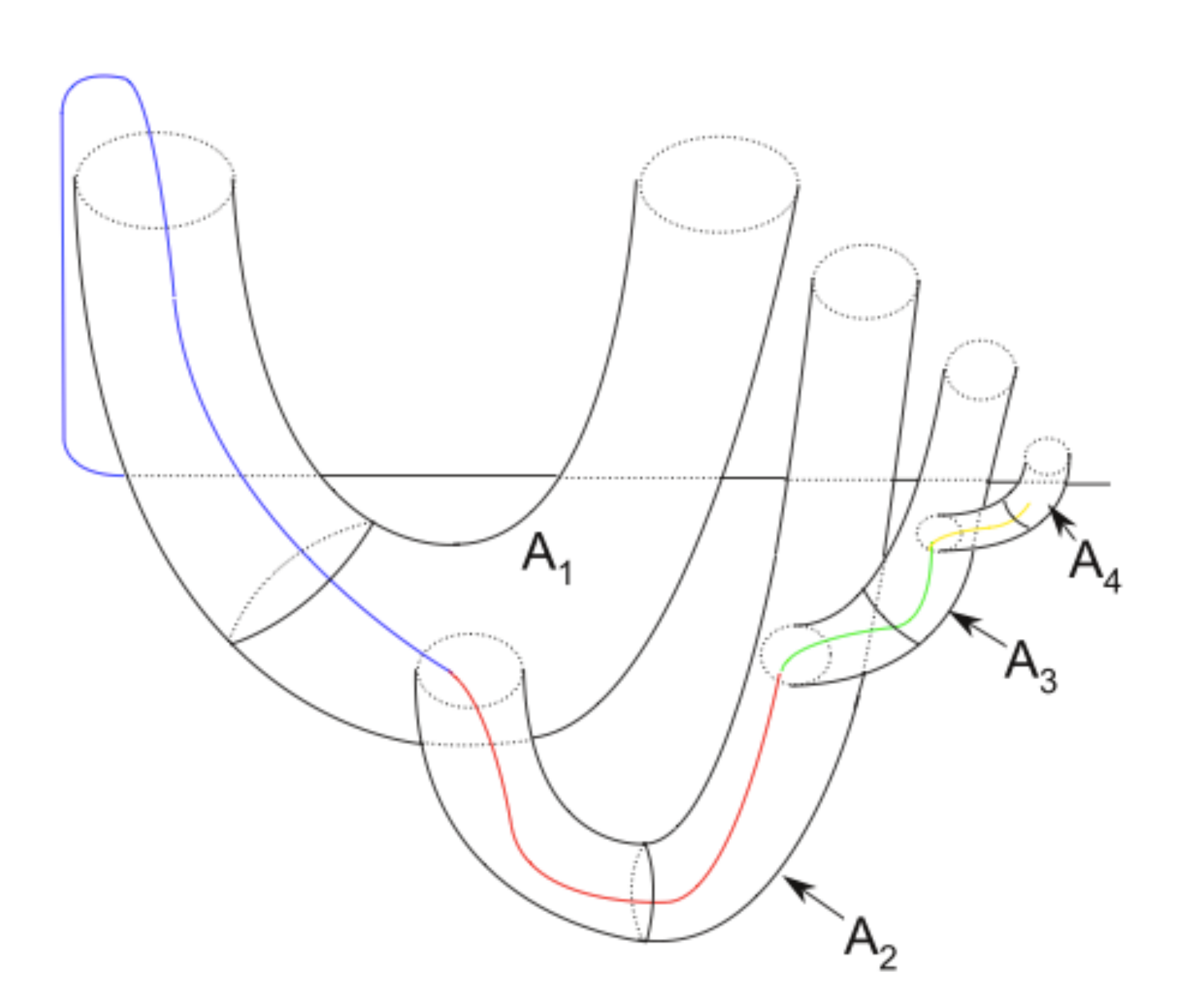}
\caption{ A picture of the first four iterations of the construction. The arc $\alpha_{i}$ terminates in the interior of the annulus $A_{i}$. }
\end{figure} 

Finally, there is a second collection of pairwise disjoint curves $\left\{d_{1},...,d_{g-2} \right\}$ such that $d_{k}$ intersects $c_{j}$ if and only if $j \leq k$.

We obtain $S_{g}$ from $S_{g-1}$ by excising a pair of open disks $D^{(g)}_{1}, D^{(g)}_{2}$ located close to the endpoints of $\alpha_{g}$, and within $N_{g-1}$, not separated by $\alpha_{g-1}$ within $N_{g-1}$, and gluing on an annulus $A_{g-1}$ along the resulting two boundary components. $c_{g-1}$ is the core curve of $A_{g-1}$, and there is a curve $d_{g-1}$ containing $\alpha_{g-1}$ as a sub-arc, which is disjoint from $d_{k}$ for all $k < g-1$, and which intersects $c_{k}$ for all $k \leq g-1$. 

We then define $\alpha_{g}$ on $S_{g}$ to be an extension of $\alpha_{g-1}$ which enters into the interior of $A_{g-1}$ and intersects $c_{g-1}$ once, and $d_{g-1}$ is a simple closed curve disjoint from $d_{k}$ for all $k \leq g-2$, and $d_{g-1}$ intersects $c_{g-1}$ once.

Note that this realization of $S_{g}$ comes equipped with a sequence of inclusions 

\[S_{1,2} \hookrightarrow S_{2,2} \hookrightarrow... \hookrightarrow S_{g-1,2} \hookrightarrow S_{g},\]

in accordance with how $A_{k}$ glues to $S_{k,2}$ to obtain $S_{k+1}$.

We are now ready to construct $\Omega(1,g)$. Define $\Omega(1,1)$ to be the single curve $\gamma_{1}$ whose isotopy class is represented by the pair of integers $(0,1)$ on $S_{1}$, and such that $\gamma_{1}$ intersects $\alpha_{1}$. $\Omega(1,g)$ is then defined by 

\[ \Omega(1,g):= \Omega(1,g-1) \cup \bigcup_{k=1}^{g-1}T_{d_{g-1}}(c_{k}),\]

where $T_{d_{g-1}}$ denotes the left Dehn-twist about $d_{g-1}$, and we think of $\Omega(1,g-1)$ as living on $S_{g}$ via the aforementioned inclusion. 

By definition, $\Omega(1,g)$ contains $g-1$ more elements than $\Omega(1,g-1)$, and therefore the sequence 
\[ \left\{|\Omega(1,g)|\right\}_{g=1}^{\infty} \]

grows quadratically as required. It remains to show that $\Omega(1,g)$ is in fact a $1$-system. For this, we must check the following four criteria:

\begin{enumerate}
\item No two elements of $\bigcup_{k=1}^{g-1}T_{d_{g-1}}(c_{k})$ are homotopic;
\item No two elements of $\bigcup_{k=1}^{g-1}T_{d_{g-1}}(c_{k})$ intersect more than once;
\item No element of $\Omega(1,g-1)$ is homotopic to any element of $\bigcup_{k=1}^{g-1}T_{d_{g-1}}(c_{k})$;
\item No element of $\Omega(1,g-1)$ intersects an element of $\bigcup_{k=1}^{g-1}T_{d_{g-1}}(c_{k})$ more than once.
\end{enumerate}

$(1)$ and $(2)$ follow from the fact that $\bigcup_{k=1}^{g-1}T_{d_{g-1}}(c_{k})$ is a homeomorphic image of a collection of pairwise disjoint, pairwise non-homotopic simple closed curves. For $(3)$, note that every element of $\bigcup_{k=1}^{g-1}T_{d_{g-1}}(c_{k})$ intersects $c_{g-1}$ essentially, but no element of $\Omega(1,g-1)$ does.

For $(4)$, we first define a many-to-one map $\Psi: \Omega(1,g) \rightarrow \left\{0,1,...,g-1\right\}$ as follows:

Orient $\alpha_{g}$ such that the endpoint it shares with $\alpha_{1}$ is the initial point of $\alpha_{g}$; this induces a compatible orientation on each sub-arc $\alpha_{h}$. Let $C_{g}:= \left\{c_{1},...,c_{g-1}\right\}$, and orient each curve in $C_{g}$ such that all intersections with $\alpha_{g}$ occur with the same orientation. Note that the index $r$ of $c_{r} \in C_{g}$ agrees with the orientation of $\alpha_{g}$, in the sense that $c_{r}$ is the $r^{th}$ element of $C_{g}$ that $\alpha_{g}$ intersects. 

Recall that $\gamma\in \Omega(1,g)\setminus \left\{\gamma_{1}\right\}$ is a Dehn twist of some curve $c_{r}$ in $C_{g}$ about $d_{k}$ for some $k$; we define $\Psi(\gamma)= r$; that is, $\Psi(\gamma)$ is the index of the curve in $C_{g}$ of which $\gamma$ is a Dehn twist. Define $\Psi(\gamma_{1}):=0$.

To show $(4)$ for a given pair of curves $\beta_{1} \in \bigcup_{k=1}^{g-1}T_{d_{g-1}}(c_{k})$ and $\beta_{2} \in \Omega(1,g-1)$, we note that up to combinatorial equivalence, there are three possible cases: either $\Psi(\beta_{1}) > \Psi(\beta_{2})$, $\Psi(\beta_{1}) \leq \Psi(\beta_{2})$, or $\Psi(\beta_{1})= \Psi(\beta_{2})$.

 Then if $\Psi(\beta_{1}) > \Psi(\beta_{2})$, $\beta_{1}$ intersects $\beta_{2}$ before arriving at $\alpha_{g}$.

Once $\beta_{2}$ is within a small neighborhood of $\alpha_{g}$, it will enter an annulus that $\beta_{1}$ does not enter, without having to cross over $\beta_{1}$. Upon returning to the corresponding element of $C_{g}$ of which $\beta_{1}$ is a Dehn-twisted image, no further intersections with $\beta_{2}$ are required because $\beta_{2}$ has already left along some earlier element of $C_{g}$ (see Figure).

\begin{figure}[H]
\centering
	\includegraphics[width=2.5in]{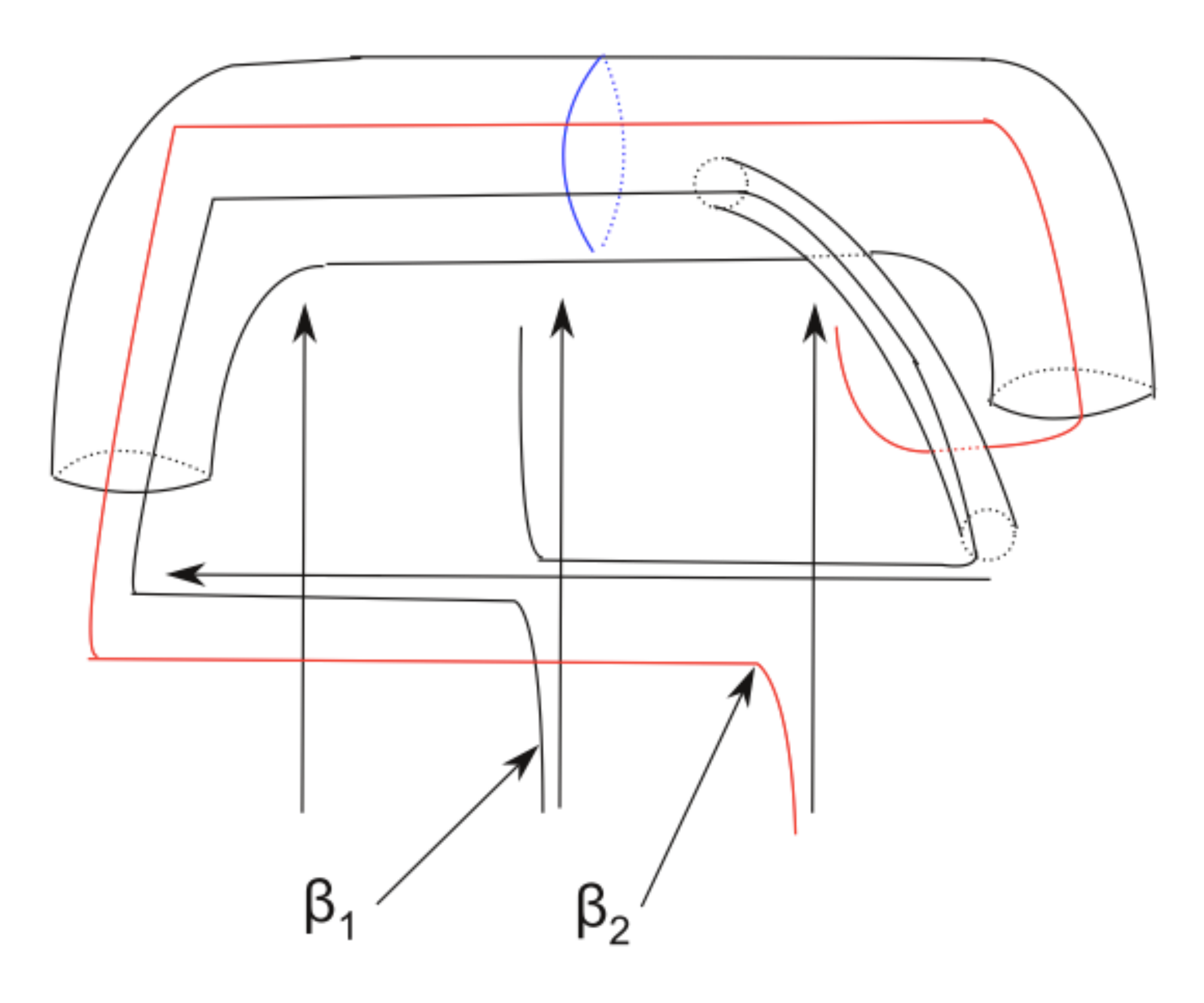}
\caption{If $\Psi(\beta_{1})>\Psi(\beta_{2})$, we can choose representatives such that $\beta_{1}$ intersects $\beta_{2}$ before arriving at $\alpha_{g}$, and never again.}
\end{figure}

\begin{figure}[H]
\centering
	\includegraphics[width=2.5in]{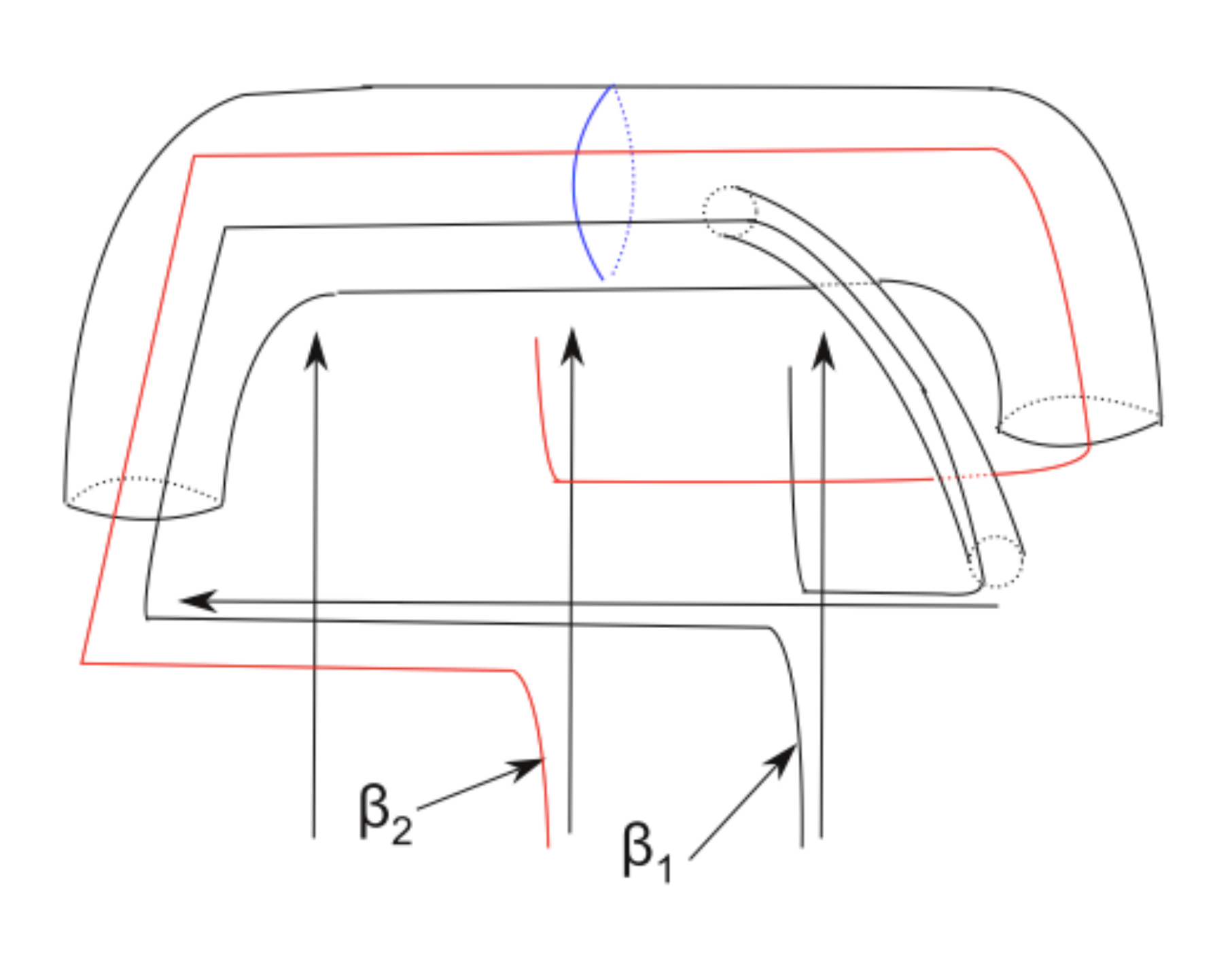}
\caption{If $\Psi(\beta_{1})<\Psi(\beta_{2})$, we can choose representatives such that $\beta_{1}$ intersects $\beta_{2}$ right before leaving $\alpha_{g}$, and never again.}
\end{figure} 

If $\Psi(\beta_{1})<\Psi(\beta_{2})$, there exists representatives of $\beta_{1}$ and $\beta_{2}$ such that $\beta_{1}$ does not intersect $\beta_{2}$ before arrival at $\alpha_{g}$; however, since $\beta_{2}$ departs from $\alpha_{g}$ later than $\beta_{1}$, $\beta_{1}$ must intersect $\beta_{2}$ once when leaving $\alpha_{g}$ (see Figure $4$).

Finally, if $\Psi(\beta_{1}) = \Psi(\beta_{2})$, there exists representatives such that $\beta_{1}$ and $\beta_{2}$ don't intersect at all within a small neighborhood of $\alpha_{g}$, but they must intersect in order to close back up.

\begin{figure}[H]
\centering
	\includegraphics[width=2.5in]{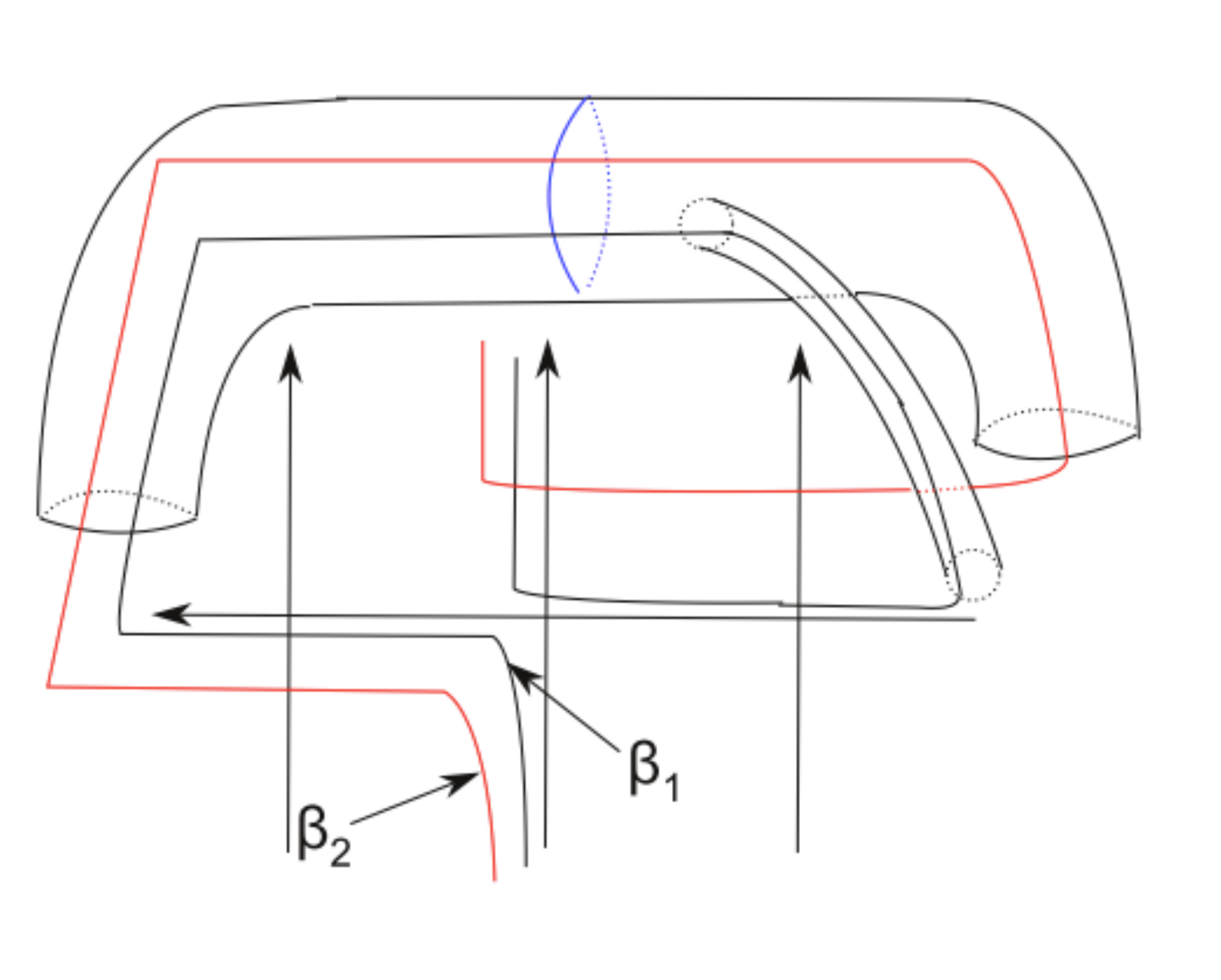}
\caption{If $\Psi(\beta_{1})=\Psi(\beta_{2})$, there exist representatives such that $\beta_{1}$ does not intersect $\beta_{2}$ near $\alpha_{g}$, but the curves must intersect in order to close up properly. }
\end{figure}

\section{Construction of $\Omega(k,g)$ for $k >1$ }

In this section, we construct $\Omega(k,g)$ for $k>1$. We will observe that our $k$-systems ``nest'' in the following sense: $\Omega(k,g)$ on $S_{g}$ is obtained from $\Omega(k,g-1)$ on $S_{g-1}$ by excising two disks, gluing an annulus on along the resulting boundary components, and adding a copy of $\Omega(k-1,g)$, ``twisted'' through the new annulus. In this way, $|\Omega(k,g)|$ will satisfy the recurrence relation 

\[ |\Omega(k,g)| = |\Omega(k,g-1)|+ |\Omega(k-1,g)|; \]

By induction, $|\Omega(k-1,g)|$ will be on the order of $g^{k}$, and therefore 

\[ |\Omega(k,g)| = \Theta\left (g^{k+1} \right). \]

However, a given pair of curves in $\Omega(k,g)$ can intersect up to $2k-1$ times, and therefore $\Omega(k,g)$ will be a $(2k-1)$-system. We therefore prove the second theorem by redefining 

\[ \Omega(k,g):= \Omega\left( \lfloor \frac{k+1}{2} \rfloor, g \right).\]

Concretely,  we define $\Omega(k,g)$ recursively by 

\[ \Omega(k,1):= \Omega(1,1), \]

and

\[ \Omega(k,g):= \Omega(k,g-1) \cup \bigcup_{\gamma \in \Omega(k-1,g)}T_{d_{g-1}}(\gamma),\]

where as in the previous section, we think of $\Omega(k,g-1)$ as living on $S_{g}$ via the sequence of inclusions described earlier. 

 As in the construction of $\Omega(1,g)$, there are $4$ requirements to varify:

\begin{enumerate}
\item No two elements of $\bigcup_{\gamma \in \Omega(k-1,g)}T_{d_{g-1}}(\gamma)$ are homotopic;
\item No two elements of $\bigcup_{\gamma \in \Omega(k-1,g)}T_{d_{g-1}}(\gamma)$ intersect more than $2k-1$ times;
\item No element of $\Omega(k,g-1)$ is homotopic to an element of $\bigcup_{\gamma \in \Omega(k-1,g)}T_{d_{g-1}}(\gamma)$;
\item No element of $\Omega(k,g-1)$ intersects an element of $\bigcup_{\gamma \in \Omega(k-1,g)}T_{d_{g-1}}(\gamma)$ more than $2k-1$ times.

\end{enumerate}

$(1)$ and $(2)$ both follow from the fact that $\bigcup_{\gamma \in \Omega(k-1,g)}T_{d_{g-1}}(\gamma)$ is a homeomorphic image of a $[2(k-1)-1]$-system, and $(3)$ follows from the fact that every element of $\bigcup_{\gamma \in \Omega(k-1,g)}T_{d_{g-1}}(\gamma)$ intersects $c_{g-1}$ essentially, but no element of $\Omega(k, g-1)$ does. 

For $(4)$, let $\beta_{1} \in \Omega(k,g-1)$, and let $\beta_{2} \in \bigcup_{\gamma \in \Omega(k-1,g)}T_{d_{g-1}}(\gamma)$; note that both $\beta_{1}$ and $\beta_{2}$ are Dehn twists of curves $\tilde{\beta}_{1}, \tilde{\beta}_{2}$, respectively, which are elements of the $[2(k-1)-1]$-system $\Omega(k-1,g)$. Concretely, $\tilde{\beta}_{1}$ is the pre-image of $\beta_{1}$ under the Dehn twist about $d_{g-1}$, and $\tilde{\beta}_{2}$ is the pre-image of $\beta_{2}$ under the Dehn twist about some $d_{k}$ for $k<g-1$. 

$i(d_{k}, d_{g-1})=0$, and 
\[ i(d_{k}, \tilde{\beta}_{1}) \leq 1 ; i(d_{g-1}, \tilde{\beta}_{2}) =1. \]

It therefore follows that 

\[ i(\beta_{1}, \beta_{2}) \leq i(\tilde{\beta}_{1}, \tilde{\beta}_{2}) +2 \leq 2k-1. \]

This completes the proof of $(4)$, and the construction of $\Omega(k,g)$.

It remains to show that 

\[ |\Omega(k,g)| \geq {g \choose 1+\lfloor\frac{k}{2} \rfloor} .\]

For this, we use the following inductive argument suggested by Kyle Luh and Daniel Montealegre:

For $k=1$, note that $\Omega(1,g)$ is obtained from $\Omega(1,g-1)$ by adding an additional $g-1$ curves, and $|\Omega(1,1)|=1$. Thus 

\[ |\Omega(1,g)|= \frac{g(g-1)}{2} > g. \]

For $k=2$, $\Omega(2,g)= \Omega(1,g)$, and  $\frac{g(g-1)}{2} = {g \choose 2}.$

Note for any $k>1$, 
\[ |\Omega(k,g)|= \sum_{i=1}^{g}|\Omega(k-1,i)|. \]

For $k$ odd, $\lfloor \frac{k}{2} \rfloor +1 = \frac{k-1}{2} +1$; by Pascal's identity, 

\[ { g \choose \frac{k-1}{2}} = {g-1 \choose \frac{k-1}{2} } + {g-1 \choose \frac{k-3}{2} } \]

\[ \leq |\Omega\left(k, g-1 \right)| + |\Omega\left(k-2 , g-1 \right)|, \]
by induction on $g$ and $k$. 

This in turn is equal to 

\[ \sum_{i=1}^{g-1} |\Omega(k-1, i)| +  |\Omega(k-2,g-1) |  \]

\[ \leq \sum_{i=1}^{g}|\Omega(k-1,i)|= |\Omega(k,g)|.\]

A similar argument holds for $k$ even; this completes the proof of Theorem $1.2$.

\end{document}